\pdfoutput=1
\RequirePackage{ifpdf}
\ifpdf 
\documentclass[pdftex]{sigma}
\else
\documentclass{sigma}
\fi

\usepackage{booktabs}

\newcommand{\N}{\mathbb{N}}
\newcommand{\Z}{\mathbb{Z}}
\newcommand{\C}{\mathbb{C}}

\newcommand{\R}{\mathbb{R}}
\newcommand{\pro}{\mathbb{P}}
\newcommand{\Po}{\mathbb{H}}

\numberwithin{equation}{section}

\newtheorem{Theorem}{Theorem}[section]
\newtheorem*{Theorem*}{Theorem}

\newtheorem{Lemma}[Theorem]{Lemma}
\newtheorem{Proposition}[Theorem]{Proposition}
 { \theoremstyle{definition}
\newtheorem{Definition}[Theorem]{Definition}

 }

\begin{document}
\allowdisplaybreaks

\newcommand{\arXivNumber}{2004.02971}

\renewcommand{\PaperNumber}{024}

\FirstPageHeading

\ShortArticleName{Accessory Parameters for Four-Punctured Spheres}

\ArticleName{Accessory Parameters for Four-Punctured Spheres}

\Author{Gabriele BOGO}

\AuthorNameForHeading{G.~Bogo}

\Address{Fachbereich Mathematik, Technische Universit\"at Darmstadt,\\ Schlossgartenstrasse 7, 64289 Darmstadt, Germany}
\Email{\href{mailto:bogo@mathematik.tu-darmstadt.de}{bogo@mathematik.tu-darmstadt.de}}

\ArticleDates{Received August 24, 2021, in final form March 22, 2022; Published online March 28, 2022}

\Abstract{We study the accessory parameter problem for four-punctured spheres from the point of view of modular forms. The value of the accessory parameter giving the uniformization is characterized as the unique zero of a system of equations. This gives an ef\-fec\-tive method to compute the uniformizing differential equation. As an application, we compute numerically and study the local expansion of the real-analytic function associating to a four-punctured sphere the value of its uniformizing parameter, and make some observations on its coefficients.}

\Keywords{accessory parameters; Fuchsian uniformization; modular forms}

\Classification{30F35; 34M03; 32G15}

\section{Introduction}
Classically, the uniformization of a genus~$g$ Riemann surface~$X$ with~$n$ punctures, $2g-2+n>0$, was related to second-order linear differential equations depending on~$3g-3+n$ parameters called \emph{accessory parameters}. Poincar\'e~\cite{Poincare} conjectured the existence of a unique choice of the accessory parameters with the following property: the ratio of linearly independent solutions of the associated differential equation lifts to a biholomorphism between the universal covering~$\widetilde{X}$ of~$X$ and the upper half-plane~$\Po$. This identification would give an explicit universal covering map for~$X$.
Despite many efforts, nobody could determine in general this choice of parameters, and the uniformization theorem was eventually proved with different techniques. The determination of the special choice of parameters is known as the~\emph{accessory parameter problem}.

Even if other approaches have proved to be better suited for the classical problem of uniformization, the accessory parameter problem is still of interest both in mathematics and physics. J.~Thompson~\cite{Thompson} discussed the algebraicity of accessory parameters of spheres with algebraic punctures in relation with Belyi's theorem; D.~Chudnovsky and G.~Chudnovsky~\cite{Chudnovsky} computed numerically the accessory parameter for genus one curves with one puncture in their numerical investigations on the Grothendiek--Katz's $p$-curvature conjecture; L.~Takhtajan and P.~Zograf~\cite{TZ1} related the accessory parameters to the Weil--Petersson metric on the Teichm\"uller space of~$n$-punctured spheres, their discoveries being stimulated by a conjecture of Polyakov's in string theory~\cite{Polyakov}. Because of the relation with Liouville theory, the computation of the accessory parameters is an active field of research in mathematical physics (see~\cite{T2} for a general introduction). Finally, a relation between local deformations and extensions of symmetric tensor representations, via accessory parameters, has been investigated in~\cite{Bogo}.

In this paper we are concerned with the simplest case of the accessory parameter problem, that of a four-punctured sphere~$X_\alpha:=\pro^1\smallsetminus\big\{\infty,1,0,\alpha^{-1}\big\}$ where~$\alpha\in\C\smallsetminus\{0,1\}$. The associated family of differential equations is of the form
\begin{equation}\label{eqn:dei}
P(t)\frac{{\rm d}^2y(t)}{{\rm d}t^2} +{P'(t)}\frac{{\rm d}y(t)}{{\rm d}t} +{(t-\rho)}y(t) =0 ,\qquad P(t) :=t(t-1)\big(t-\alpha^{-1}\big),
\end{equation}
where~$\rho\in\C$ is the~\emph{accessory parameter}. We call the unique choice~$\rho_{\rm F}$ of the accessory parameter inducing the identification~$\widetilde{X}_\alpha\to\Po$ the~\emph{Fuchsian value}.

Apart from very special cases, e.g., when the differential equation~\eqref{eqn:dei} is a Picard--Fuchs differential equation~\cite{BouwMoeller, ZagApery}, it is still not known how to determine the Fuchsian parameter even in this simplest case. Several papers in the literature deal with the numerical computation of the Fuchsian parameter for four-punctured spheres or, equivalently, elliptic curves with one puncture. Other to the above mentioned work of Chudnovsky and Chudnovsky, one should mention the work of L.~Keen, H.~Rauch, and A.~Vasquez~\cite{KRV}, and J.~Hoffman's Ph.D.~Thesis~\cite{Hoffman}.
These works are based on the observation that the monodromy group of the uniformizing differential equation, which is the Deck group of the universal covering~$\Po\to X$, is a discrete subgroup of~${\rm SL}_2(\R)$. To require the monodromy of~\eqref{eqn:dei} to have real coefficients imposes constraints on the choice of the accessory parameters that can be used to numerically compute them. However, since~\eqref{eqn:dei} can have a monodromy group with real coefficients without being the uniformizing differential equation (in fact, this happens for a discrete set of accessory parameters), a further analysis to determine the Fuchsian one is needed. A new idea, based on the theory of Painlev\'e~VI equation and isomonodromy deformations, leading to the numerical computation of the Fuchsian parameter has recently appeared in~\cite{Anselmo}.

In this note a different approach, based on the modularity of the solution of the uniformizing differential equation, is described. As an application, we get an efficient way to compute numerically the value of the Fuchsian parameter of a given four-punctured sphere~$X_\alpha$ in terms of~$\alpha$.
The main result can be stated as follows
\begin{theorem*}[Theorem~\ref{thm:main}]
The Fuchsian value for the punctured sphere~$X_\alpha$ is the unique zero of a system of infinitely many equations constructed from the differential equation~\eqref{eqn:dei}.
\end{theorem*}
What makes the accessory parameter problem hard is that the dependence of the uniformization data (monodromy, covering map) on the location of the punctures is quite obscure. Theorem~\ref{thm:main} shows that in the case of four-punctured spheres~$X_\alpha$ it is possible to construct a system of equations solved by the Fuchsian parameter using only basic properties of~$X_\alpha$, namely the existence of non-trivial automorphisms.
Nehari~\cite{Nehari}, using different ideas, also characterized the Fuchsian parameter as a zero of a system of infinitely many equations in the case all the punctures lie on the real line.

We describe the main ideas of the proof of Theorem~\ref{thm:main}.
\begin{enumerate}\itemsep=0pt
\item For every choice of~$\alpha\in\C\smallsetminus\{0,1\}$, the punctured sphere~$X_\alpha\simeq\Po/\Gamma_\alpha$ has a Klein group of automorphisms permuting the punctures. The fixed points of these automorphisms are in correspondence with cusp representatives of the uniformizing group~$\Gamma_\alpha$. It turns out that a set of generators of~$\Gamma_\alpha$ can be described purely in terms of these cusp representatives and then, via the covering map, in terms of fixed points of~$X_\alpha$. This is discussed in Section~\ref{sec:four}. From the point of view of the differential equation~\eqref{eqn:dei} we have the following description. If~$\rho=\rho_{\rm F}$ is the Fuchsian parameter, the ratio~$\eta_{\rho_{\rm F}}$ of independent solutions is an inverse of the covering map; the images of the fixed points of~$X_\alpha$ via~$\eta_{\rho_{\rm F}}$ can be used to construct the cusp representatives of~$\Gamma_\alpha$ and finally the uniformizing group itself. It follows that the group constructed in this way is the monodromy group of~\eqref{eqn:dei} for~$\rho=\rho_{\rm F}$.

\item For a generic choice of the accessory parameter~$\rho$, the ratio~$\eta_\rho$ of independent solutions of~\eqref{eqn:dei} is not an injective map. However, there is an open set~$\mathcal{B}$ of parameters such that~$\eta_\rho$ is injective if~$\rho\in\mathcal{B}$ (such~$\eta_\rho$ gives a quasi-Fuchsian uniformization of~$X_\alpha$). Of course \mbox{$\rho_{\rm F}\in\mathcal{B}$}. The idea is to mimic the construction of the previous paragraph for the accessory parameters in~$\mathcal{B}$. More precisely, we attach a Fuchsian group~$\Gamma(\rho)$ to every~$\rho$ in an open subset of~$\mathcal{B}$ by defining real numbers~$c_1(\rho)$,~$c_2(\rho)$ (``potential'' cusp representatives) from the images of the fixed points of~$X_\alpha$ via~$\eta_\rho$. The group~$\Gamma(\rho)$ is constructed from~$c_1(\rho)$, $c_2(\rho)$ by using Poincar\'e's theorem. This is discussed in Section~\ref{sec:find}. We remark that contrary to the case~$\rho=\rho_{\rm F}$, the group~$\Gamma(\rho)$ is not the monodromy group of~\eqref{eqn:dei} which, for~$\rho_{\rm F}\neq\rho\in\mathcal{B}$ is not a Fuchsian group.

\item When~$\rho=\rho_{\rm F}$, a holomorphic solution of~\eqref{eqn:dei} lifts to a holomorphic modular form~$f(\tau)$ for the uniformizing group~$\Gamma_\alpha$; its~$q$-expansion is easily computed from the solutions of~\eqref{eqn:dei}. Similarly, for any~$\rho$ we can construct a~$Q$-expansion~$f_\rho(Q)$ whose coefficients depend on~$\rho$. For~$\rho$ in a subset of~$\mathcal{B}$ we can test the modularity of~$f_\rho(Q)$ with respect to the Fuchsian group~$\Gamma(\rho)$ constructed in the previous step. It turns out that~$f_\rho(Q)$ is modular for~$\Gamma(\rho)$ if and only if~$\rho=\rho_{\rm F}$ (Section~\ref{sec:findfin}). The equations describing the modular transformations of~$f_\rho(Q)$ with respect to the generators of~$\Gamma(\rho)$ give the system in Theorem~\ref{thm:main}.
\end{enumerate}

As mentioned above, this construction gives an efficient method to compute numerically the Fuchsian parameter by approximating a solution of the system of equations in Theorem~\ref{thm:main}. Section~\ref{sec:xpl} presents an application of this method to the study of the analytic properties of the Fuchsian accessory parameter function. This map associates to the four-punctured sphere~$X_\alpha$ its Fuchsian value~$\rho_{\rm F}(X_\alpha)$; we can see this as a map~$\rho_{\rm F}\colon\C\smallsetminus\{0,1\}\to\C$. It is known that this map is real-analytic and not holomorphic. By using the method presented above we computed the coefficients of its local expansion for different values of~$\alpha$. An interesting phenomenon we can observe from the numerical data (tables at page~\pageref{table1}) concerns the size of the coefficients of this expansion. It appears that the holomorphic part of the Fuchsian parameter function~$\rho_{\rm F}$ is much larger than the rest. This suggests that the function~$\rho_{\rm F}$ may have nice analytic properties, for instance be a quasiregular map.

\section{Uniformization, modular forms, and differential equations}\label{sec:pre}
\subsection{Uniformization and differential equations}
We recall the classical theory in the case of hyperbolic Riemann surfaces of genus zero. A good reference for the general theory is the book~\cite{SGervais}.
Let $X:=\pro^1\smallsetminus\{\alpha_1,\alpha_2,\dots,\alpha_{n-2},\alpha_{n-1}=0,\alpha_n=\infty\}$, where $\alpha_i\in\C\smallsetminus\{0\}$ and $\alpha_i\neq \alpha_j$ if~$i\neq j$, be an~$n$-punctured sphere.
Consider a second-order linear differential equation on~$X$ with holomorphic coefficients:
\[
\frac{{\rm d}^2y(t)}{{\rm d}t^2}+p(t)\frac{{\rm d}y(t)}{{\rm d}t}+q(t)y(t)=0.
\]
Let~$y_1(t)$ and $y_2(t)$ be linearly independent solutions. The ratio $\eta(t):=y_2(t)/y_1(t)$ can be analytically continued to the Riemann surface~$X$ and induces a non-constant function \mbox{$\widetilde{\eta}(t)\colon\widetilde{X}\to\C$} on the universal covering $\widetilde{X}$ of~$X$. It is easy to verify that~$\widetilde{\eta}$ is a local biholomorphism. Conversely, every local biholomorphism $\widetilde{X}\to\C$ arises in this way. In particular, every global biholomorphism between $\widetilde{X}$ and a subdomain of~$\C$, if any, arises from the ratio of linearly independent solutions of differential equations of the form \cite{SGervais, Ford}
\begin{equation}\label{eqn:unifschw}
\frac{{\rm d}^2y(t)}{{\rm d}t^2} + \left(\frac{1}{4}\sum_{j=1}^{n-1}{\frac{1}{(t- \alpha_j)^2}}+\frac{1}{2}\sum_{j=1}^{n-1}{\frac{m_j}{(t-\alpha_j)}}\right)y(t)=0 ,
\end{equation}
where $m_0,\dots,m_{n-1}$ are complex parameters, called \emph{accessory parameters}, subject to the following relations\footnote{In the literature often appears another parameter $m_n$ associated to the puncture at $\infty;$ it is defined from the asymptotic expansion of the rational function in \eqref{eqn:unifschw} as $t\to\infty$. It turns out that $m_n$ can be expressed in terms of $m_1,\dots,m_{n-1}$ and of the punctures $\alpha_1,\dots,\alpha_{n-1}$ as $m_n=\sum_{j=1}^{n-1}{\alpha_j(1+m_j\alpha_j)}$.}
\begin{equation} \label{eqn:rel}
\sum_{j=1}^{n-1}{m_j}=0,\qquad\sum_{j=1}^{n-1}{\alpha_jm_j}=1-\frac{n}{2}.
\end{equation}
More precisely, for certain choices of the accessory parameters~$m_1,\dots,m_{n-1}$, the ratio of linearly independent solutions induce a biholomorphic map~$\widetilde{\eta}\colon\widetilde{X}\to \widetilde{\eta}\big(\widetilde{X}\big)\subset\C$.

The name ``accessory parameters'' is due to the fact that the choice of $m_1,\dots,m_{n-1}$ does not affect the local behaviour of solutions of~\eqref{eqn:unifschw} near the singular points (but of course influences the global behaviour of the solutions).

Nowadays it is well known that the space of accessory parameters inducing a biholomorphism is a non-empty open connected set~$\mathcal{B}$ in~$\C^{n-3}$ called the \emph{Bers slice}~\cite{Bers}. The image of the map~$\widetilde{\eta}\colon\widetilde{X}\to\C$ is in general a \emph{quasidisk}, i.e., the image of a disk under a quasiconformal transformation, with a nowhere-smooth boundary of Hausdorff dimension~$>1$. However, for a special choice of the accessory parameters the universal covering~$\widetilde{X}$ is identified, via~$\widetilde{\eta}$, with the upper half-plane~$\Po$. It is known since Poincar\'e~\cite{Poincare} that this choice is unique; we call the unique value of the accessory parameters giving the above identification the \emph{Fuchsian value}, and the corresponding differential equation the~\emph{uniformizing differential equation}. The monodromy group of the uniformizing differential equation is the Deck group of transformations of the covering; it follows that it is conjugated to a discrete subgroup of~${\rm SL}_2(\R)$ (but the converse is not true: there exists infinitely many choices of the accessory parameters such that the monodromy group is discrete in~${\rm SL}_2(\R)$, but the ratio of solutions does not induce a biholomorphic map on the universal covering~\cite{Goldman}.)

The \emph{accessory parameter problem} consists in finding the Fuchsian value for a given punctured sphere~$X$. This problem turned out to be very hard and only partial or numerical solutions for spheres with a low number of punctures (in fact, only four punctures) have been found. It is worth noting that even the existence of the Fuchsian value has never been proved directly, i.e., without referring to the uniformization theorem; the only exception is the case of four-punctured spheres with real punctures, which was solved by V.~Smirnov~\cite{Smirnov}.

\subsection{Modular forms and differential equations}\label{sec:funcor}

Let~$\Gamma\subset{\rm SL}_2(\R)$ be a cofinite discrete group, let~$t\colon\Po/\Gamma\to\C$ be a modular function, and let~$f\in M_k(\Gamma)$ be a modular form of weight~$k$ on~$\Gamma$. If we express locally~$f$ as a function of~$t$, i.e., $\phi(t(\tau))=f(\tau)$, then the function~$\phi(t)$ satisfies a linear differential equation of order~$k+1$ with algebraic coefficients. Similarly, a~$k$-th root of~$f$, if expressed as a function of~$t$, satisfies a linear differential equation of order~$2$ with algebraic coefficients; a local basis of solutions is given by~$\big\{\phi(t),\hat{\phi}(t)\big\}$ where~$\phi(t(\tau))=f^{1/k}(\tau)$, $\hat{\phi}(t(\tau))=\tau f^{1/k}(\tau)$ (see Chapter~5 of the first part of \cite{123} for details).

Now let~$\Gamma$ be of genus zero and torsion free, and let~$t$ be a Hauptmodul, i.e., a modular function that extends to an isomorphism between the compactification of~$\Po/\Gamma$ and~$\pro^1(\C)$. In this setting the linear differential equation satisfied by~$f^{1/k}$ is defined on the punctured sphere~$t(\Po/\Gamma)$ and its coefficients are rational functions of~$t$. Since the ratio of the independent solutions~$\phi(t)$,~$\hat{\phi}(t)$ lifts to the coordinate~$\tau\in\Po$, we see that the differential equation satisfied by a~$k$-th root of (every)~$f\in M_k(\Gamma)$ with respect to~$t$ is the uniformizing differential equation of~$t(\Po/\Gamma)$ in the sense of the previous section. We can then reformulate the accessory parameter problem as follows.
\begin{Proposition}\label{thm:refor}
The Fuchsian value is the unique choice of accessory parameters such that the holomorphic solution of the associated differential equation lifts to a~$k$-th root of a modular form~$f\in M_k(\Gamma)$ with respect to the monodromy group~$\Gamma\subset{\rm SL}_2(\R)$.
\end{Proposition}

If the Hauptmodul~$t$ is fixed, different choices of~$f\in M_*(\Gamma)$ yield different differential equations; however, the ratio of independent solutions always lift to~$\tau\in\Po$, that is, differential equations associated to different choices of~$f$ are projectively equivalent. In particular, the equation in~\eqref{eqn:unifschw} correspond to the choice of the meromorphic modular form~$f={\rm d}t/{\rm d}\tau$. A way to see this is to describe the coefficients of the differential equation in terms of \emph{Rankin--Cohen brackets}: if~$f\in M_k(\Gamma)$ and~$g\in M_l(\Gamma)$ these are defined as follows
\[
[f,g]_1:=kfg'-lgf' ,\qquad [f,g]_2:=\frac{k(k+1)}{2}fg''-2(k+1)(l+1)f'g'+\frac{l(l+1)}{2}gf'',
\]
where~$'=(2\pi {\rm i})^{-1}{\rm d}/{\rm d}\tau$. The brackets~$[f,g]_1$ and~$[f,g]_2$ are modular forms of weight~$k+l+2$ and~$k+l+4$ respectively. If we set~$g=t'$, which is a meromorphic modular form of weight~$2$, the quotients
\[
A(t) :=\frac{[f,t']}{kft'^2},\qquad B(t) :=-\frac{[f,f]_2}{k^2(k+1)f^2t'^2}
\]
are modular forms of weight zero, so in particular they are rational functions of the Hauptmodul~$t$. It is easy to verify that the differential equation satisfied by~$\phi(t)=f^{1/k}(\tau)$ is given by~${\rm d}^2\phi(t)/{\rm d}t^2+A(t){\rm d}\phi(t)/{\rm d}t+B(t)\phi(t)=0$.
In the case also $f=t'$ a simple computation reveals that~$A(t)=0$ and~$B(t)$ is the Schwarzian derivative~$B(t)=\{\tau,t\}/2$, where~$\tau\in\Po$. We can then recall the classical identity (see for example the first section of~\cite{TZ1})
\[
\{\tau,t\}=\frac{1}{2}\sum_{j=1}^{n-1}{\frac{1}{(t- \alpha_j)^2}}+\sum_{j=1}^{n-1}{\frac{m_j}{(t-\alpha_j)}}
\]
to conclude that the differential equation satisfied by~$\phi(t)=\sqrt{t'}$ is precisely~\eqref{eqn:unifschw}.

In Appendix~\ref{appendixA} we compute the differential equation associated to (the square root of) a~special choice of~$f\in M_2(\Gamma)$. In the case~$n=4$ this reduces to the well-known Heun equation.

\section{Four-punctured spheres}\label{sec:four}
In this section we show how the generators of a torsion-free Fuchsian group with four cusps and the automorphisms of the corresponding punctured sphere are related.

\begin{figure}[t]\centering

\includegraphics[width=0.4\textwidth]{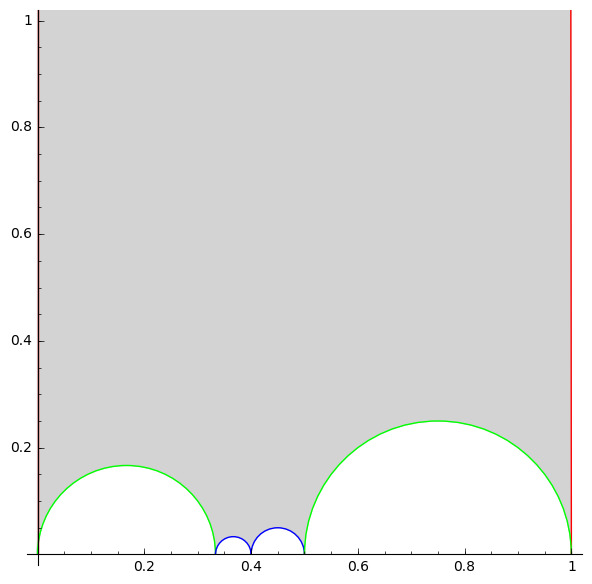}
\caption{Fundamental domain of~$\Gamma_1(5)$. The cusp representatives are~$0,1/3,2/5$ and $\infty$.}\label{fig:fd}
\end{figure}

 Let $\alpha\neq 0,1$ be a complex number and consider the four-punctured sphere~$X_\alpha:=\pro^1\smallsetminus\big\{\infty,1,0,\alpha^{-1}\big\}$.
We are going to choose a uniformizing group~$\Gamma_\alpha$ and a Hauptmodul~$t$ for~$X_\alpha$. A~priori they are not uniquely defined: the uniformizing group is determined only up to conjugacy in~${\rm SL}_2(\R)$, and the composition of a given Hauptmodul with any automorphism of~$X$ still yields a Hauptmodul.
In any case, the group $\Gamma_\alpha$ is torsion free and has four non-equivalent cusps;
we denote the equivalence classes of cusps by~$[c_1]$, $[c_2]$, $[c_3]$, $[c_4]$ (later we will fix $[c_3]=[\infty]$ an $c_4=[0])$.
The cusps are in bijection, via~$t$, with the punctures of $X_\alpha$.
A picture of a fundamental domain for the action of~$\Gamma_\alpha$ on $\Po$ in the special case~$\Gamma$ is the congruence subgroup~$\Gamma_1(5)$ is given in Figure~\ref{fig:fd}.
The next lemma describes the normalization of~$\Gamma_\alpha$ and~$t$ we choose.

\begin{Lemma}\label{thm:norm}
Let $X_\alpha=\pro^1 \smallsetminus \big\{\infty,1,0,\alpha^{-1}\big\}$, $\alpha\in\C\smallsetminus\{0,1\}$.
There exists a pair~$(\Gamma_\alpha,t)$ with~$t\colon\overline{\Po/\Gamma_\alpha}\overset{\sim}{\rightarrow}\overline{X_\alpha}=\pro^1(\C)$ and such that
\begin{enumerate}\itemsep=0pt
\item[$1)$] the group~$\Gamma_\alpha$ has inequivalent cusps~$[\infty]$, $[0]$ and the stabilizer~${\Gamma_{\alpha}}_\infty$ of~$\infty$ in~$\Gamma_\alpha$ is generated by
\[
T=\begin{pmatrix}
1 & 1\\
0 & 1
\end{pmatrix}\in\Gamma_\alpha ;
\]
\item[$2)$]
the values of $t$ at the inequivalent cusps $\infty,0$ are~$t(\infty)=0$ and $t(0)=\alpha^{-1}$.
\end{enumerate}
These choices uniquely determine $\Gamma$ and $t$.
\end{Lemma}
\begin{proof}
Let $\hat{\Gamma}$ and $\hat{t}$ be such that $X_\alpha=\hat{t}\bigl(\Po/\hat{\Gamma}\bigl)$. If $\hat{\Gamma}$ and $\hat{t}$ are as in the statement we are done. If not, compose~$\hat{t}$ with an automorphism~$\phi$ of~$X_\alpha$ in such a way that $\hat{t}_1:=\phi\circ\hat{t}$ maps the cusp $\infty$ to the puncture $0$ (such automorphism of $X_\alpha$ always exists, see the next section).
For nonzero real numbers~$a$, $b$ consider the matrix~$\sigma:=\bigl(\begin{smallmatrix}
a & b\\
0 & 1/a
\end{smallmatrix}\bigr)$ and define
\[
\Gamma(a,b):=\sigma\hat{\Gamma}\sigma^{-1},\qquad t(\tau):=\hat{t}_1\big(\sigma^{-1}\tau\big) .
\]
The map~$t$ is by construction a Hauptmodul for~$\Gamma(a,b)$; since~$\sigma\infty=\infty$ and~$\hat{t}_1(\infty)=0$ we also have~$t(\infty)=0$.
The conjugation of~$\hat{\Gamma}$ by~$\sigma$ amounts to determine the coordinate~$\tau$ on~$\Po$ given by the uniformizing differential equation up to a linear map~$\tau\mapsto a\tau+b$, where~$a,b\in\R$. To fix it uniquely we only need to choose~$a$ and~$b$.
A simple computation shows that the generator of~$\Gamma(a,b)_\infty$, the stabilizer of~$\infty$ in~$\Gamma(a,b)$, only depends on the choice of~$a;$ we can choose~$a=\bar{a}$ such that~$\Gamma(\bar{a},b)_\infty=\langle{T}\rangle$.
Finally, we see that $t(0)=\hat{t}_1\big(\sigma^{-1}0\big)=\hat{t}_1(b)$, and we can pick any $\bar{b}\in\R$ such that~$\hat{t}_1(\bar{b})=\alpha^{-1}$. Then $\big(\Gamma_\alpha:=\Gamma(\bar{a},\bar{b}),t\big)$ is the desired pair.
\end{proof}

In the following, $\Gamma_\alpha$ and $t$ will always be normalized as in the above lemma. In this case, the Fourier expansion of $t$ at $\infty$ starts
\begin{equation}\label{eqn:expt}
t=rq+\cdots,\qquad q={\rm e}^{2\pi {\rm i}\tau},\quad \tau\in\Po,
\end{equation}
for some~$r\in\C\smallsetminus\{0\}$.

\subsection{Generators of the uniformizing group}\label{sec:unifgrp}
The goal of this section is to write a set of parabolic generators of a torsion-free genus zero Fuchsian group $\Gamma$ with four cusps only in terms of cusp representatives. By a set of parabolic generators we mean a set of matrices~$\{M_1,\dots,M_4\}$ with~$\operatorname{Tr}^2(M_i)=4$ that generate~$\Gamma$ with the relation~$\prod_{i=1}^4{M_i}={\rm Id}$.

It follows from the existence of non-trivial automorphisms of~$\Po/\Gamma$ that the cusps of~$\Gamma$ are all regular or irregular.
In the next lemma (and in the rest of the paper) we assume that the cusps are regular; the case of irregular cusps can be handled analogously.

\begin{Lemma}\label{thm:rep}
Let~$\Gamma$ be a torsion free Fuchsian group of genus zero with four cusps. Assume that~$T=\bigr(\begin{smallmatrix}1 & 1 \\0 & 1\end{smallmatrix}\bigr)\in\Gamma$, and let~$0<c_1<c_2$ be representatives of the non-equivalent finite cusps. Then~$\Gamma=\big\langle T,S_0,S_{c_1},S_{c_2}\,|\,S_{c_2}S_{c_1}S_0T^{-1}={\rm Id}\big\rangle$ where
\begin{gather*}
S_0 =\begin{pmatrix}
1 & 0 \\
D_0 & 1
\end{pmatrix},\qquad \!
S_{c_1} =\begin{pmatrix}
1+c_1D_{c_1} & -c_1^2D_{c_1} \\
D_{c_1} & 1-c_1D_{c_1}
\end{pmatrix},\qquad\!
S_{c_2} =\begin{pmatrix}
1+c_2D_{c_2} & -c_2^2D_{c_2} \\
D_{c_2} & 1-c_2D_{c_2}
\end{pmatrix},
\end{gather*}
and the constants $D_0$, $D_{c_1}$, $D_{c_2}$ are given by
\begin{gather}\label{eqn:Di}
D_0=\frac{1}{c_1(1-c_2)},\qquad D_{c_1}=\frac{1}{c_1(c_2-c_1)},\qquad D_{c_2} =\frac{1}{(c_2-c_1)(1-c_2)} .
\end{gather}
\end{Lemma}
\begin{proof}It is well known that~$\Gamma$ is generated by the stabilizers of its cusps, and that the stabilizer of the finite (regular) cusp~$c_i$ is of the form
\begin{equation}\label{eqn:stab}
S_{c_i}=\begin{pmatrix}
1+c_iD_{c_i} & -c_i^2D_{c_i} \\
D_{c_i} & 1-c_iD_{c_i}
\end{pmatrix}
\end{equation}
for some positive $D_{c_i}\in\R$. The only thing to prove are the formulae in~\eqref{eqn:Di}.

The choice of cusp representatives in the statement fixes a fundamental domain $\mathcal{F}$ for the action of $\Gamma$.
It is well known that a free generating set for~$\Gamma$ is given by the M\"obius transformations which pairs the boundary geodesics of~$\mathcal{F}$.
Among these transformations there is one that fixes one of the finite cusp representatives (see for instance Figure~\ref{fig:fd}, where this cusp representative is~$2/5$).
In our case, the fixed cusp representative is~$c_2$, since we have set $0<c_1<c_2$. Call~$S_{c_2}$ the transformation that fixes~$c_2$ and pairs the relative boundary geodesics.

The transformation~$S_{c_2}$ also exchanges~$c_1$ with its equivalent $c_1'>c_2$. There is a transformation that exchanges~$c_1$ with~$c_1'$, and sends~$0$ to~$1$; call it~$P_{0,c_1}$. It follows that
\[
S_{c_1}:=S_{c_2}^{-1}P_{0,c_1}
\]
fixes $c_1$. In the same way $S_{0}:=P_{0,c_1}^{-1}T$ fixes $0$.
The matrices~$S_*$, $*=0,c_1,c_2$ generate the stabilizer of the cusp~$*$ and satisfy the parabolic relation~$S_{c_2}S_{c_1}S_0T^{-1}={\rm Id}$.
It follows that every $S_*$, is of the form~\eqref{eqn:stab} and then one can compute the real numbers $D_*$, $*=0,c_1,c_2$, by solving the system given by the parabolic relation
\[
\begin{pmatrix}
1+c_2D_{c_2} & -c_2^2D_{c_2} \\
D_{c_2} & 1-c_2D_{c_2}
\end{pmatrix}
\begin{pmatrix}
1+c_1D_{c_1} & -c_1^2D_{c_1} \\
D_{c_1} & 1-c_1D_{c_1}
\end{pmatrix}
\begin{pmatrix}
1 & 0 \\
D_0 & 1
\end{pmatrix} =
\begin{pmatrix}
1 & 1 \\
0 & 1
\end{pmatrix}.
\]
The formulae in~\eqref{eqn:Di} follow after an easy computation.
\end{proof}

\subsection{Automorphisms of four-punctured spheres and cusp representatives}\label{sec:cusprep}

For every choice of $\alpha\neq 0,1$, the surface $X_\alpha$ admits a Klein four-group of automorphisms $\operatorname{Aut}_0(X_\alpha)$ generated by any two of the involutions
\begin{equation}\label{eqn:autom}
\phi_0\colon\ t\mapsto\frac{1-\alpha t}{\alpha(1-t)},\qquad\phi_1\colon\ t\mapsto\frac{t-1}{\alpha{t}-1},\qquad \phi_2\colon\ t\mapsto\frac{1}{\alpha t},
\end{equation}
where $\phi_0=\phi_1\circ\phi_2$. In general $\operatorname{Aut}(X_\alpha)=\operatorname{Aut}_0(X_\alpha)$, but for exceptional choices of $\alpha$, the automorphism group of $X_\alpha$ is larger. If $\alpha=-1,1/2,2$ then $\operatorname{Aut}(X_\alpha)$ has order~$8;$ if $\alpha=1/2\pm {\rm i}\sqrt{3}/2$ then $\operatorname{Aut}(X_\alpha)$ has order~$12$. In these exceptional cases, the Fuchsian parameter and the uniformization of~$X_\alpha$ can be easily computed (see~\cite{Hempel}).

Let $t\colon\Po/\Gamma\to X_{\alpha}$ be normalized as in Lemma~\ref{thm:norm}. Every automorphism $\phi\in\operatorname{Aut}(X_\alpha)$ lifts to an automorphism $\tilde{\phi}$ of the universal covering $\Po$.
Every such automorphism~$\tilde{\phi}$ can be represented by a matrix $W_\phi\in{\rm SL}_2(\R)$ and belongs to the normalizer $N(\Gamma)$ of $\Gamma$. This, together with Lemma~\ref{thm:norm}, implies that $W_\phi TW_\phi^{-1}\in\Gamma$.
In particular, if $\phi\in\operatorname{Aut}_0(X_\alpha)$ and~$W_\phi$ sends the cusp~$[c]$ to the cusp~$\infty$, the element $S_c:=W_\phi TW_\phi^{-1}$ sends the cusp~$[c]$ to itself, so it belongs to the stabilizer~$\Gamma_c$. Actually, more is true.
\begin{Lemma}\label{thm:stab}
Let~$W_\phi$ be an involution of~$\Po$ obtained by lifting~$\phi\in\operatorname{Aut}_0(X_\alpha)$ and such that $W_\phi(c)=\infty$ for~$[c]$ a finite cusp of~$\Gamma$.
Then the transformation $S_c=W_\phi TW_\phi^{-1}$ generates the stabilizer $\Gamma_c$ of~$[c]$ in~$\Gamma$.
\end{Lemma}

From Lemma~\ref{thm:stab} and~\eqref{eqn:stab} it follows that~$W_\phi$ is of the form
\begin{equation*}
W_\phi = \sqrt{D_c}
\begin{pmatrix}
c & \dfrac{-1-c^2D_c}{D_c} \\
1 & -c
\end{pmatrix}
\end{equation*}
for some representative~$c$ of~$[c]$ and some positive constant~$D_c$.
In other words, we can describe every lift~$W_\phi$ of $\phi\in\operatorname{Aut}_0(X_\alpha)$ in terms of the cusp representatives~$c_1$, $c_2$ and the positive real constants $D_0$, $D_{c_1}$, $D_{c_2}$ in~\eqref{eqn:Di}.
The transformation~$W_\phi$ has a unique fixed point in~$\Po$ given by
\begin{equation}\label{eqn:fph}
\tau_\phi=c+{\rm i}/\sqrt{D_c} ;
\end{equation}
its image via~$t$ is a fixed point of the involution~$\phi\in\operatorname{Aut}_0(X_\alpha)$. The next lemma establish which fixed points of the automorphisms~$\phi_j\in\operatorname{Aut}_0(X_\alpha)$, $j=0,1,2$, are images of the fixed points~$\tau_j$ of~$W_{\phi_j}$.
In the next lemma, and in the rest of the paper, we assume that~$\alpha$ satisfies
\begin{equation}\label{eqn:norma}
\bigr|\alpha^{-1}\bigl|\le 1\qquad\text{and}\qquad\operatorname{Re}\bigl(\alpha^{-1}\bigr),\operatorname{Im}\bigl(\alpha^{-1}\bigr)\ge0 .
\end{equation}
One can always reduce to this case via conformal transformations and complex conjugation. The behavior of the accessory parameters with respect to these transformations is known (see~\cite{Hempel}).

\begin{Lemma}\label{thm:fixed}
Let~$\alpha^{-1}$ be as in~\eqref{eqn:norma} and let~$t\colon\Po/\Gamma\to X_\alpha$ be as in Lemma~{\rm \ref{thm:norm}}.
Let~$W_{\phi_j}$ be the lifting of~$\phi_j\in\operatorname{Aut}_0(X_\alpha)$, $j=0,1,2$. Then the image~$z_j\in X_\alpha$ of the unique fixed point~$\tau_j$ of~$W_{\phi_j}$ on~$\Po$ is given, in terms of~$\alpha$, by
\[
z_0 =\frac{\alpha-\sqrt{\alpha(\alpha-1)}}{\alpha} ,\qquad z_1=\frac{1+\sqrt{1-\alpha}}{\alpha} ,\qquad z_2 =\frac{1}{\sqrt{\alpha}} .
\]
\end{Lemma}
\begin{proof}
The fixed points of $\phi_0$ are the solutions of
\begin{equation}
\label{eqn:root}
\alpha{t}^2-2\alpha{t}+1=0.
\end{equation}
Since~$\phi_0(0)=\alpha^{-1}$ and~$t(\infty)=0$, $t(0)=\alpha^{-1}$, the lift~$W_{\phi_0}$ sends the cusp~$\infty$ to~$0$. It follows that the fixed point of~$W_{\phi_0}$ on~$\Po$ is $\tau_{0}={\rm i}/\sqrt{D_0}$. As~$\tau_{0}$ lies on the imaginary axis, its image on~$X_\alpha$ belongs to the geodesic (determined by~$t$) joining the punctures~$t(\infty)=0$ and~$t(0)=\alpha^{-1}$. Looking at the two roots of \eqref{eqn:root} and considering the constraints~\eqref{eqn:norma} on~$\alpha^{-1}$ it follows that
\[
t(\tau_{0})= z_0 =1-\frac{\sqrt{\alpha(\alpha-1)}}{\alpha}.
\]

Now consider the involution $\phi_1\in\operatorname{Aut}_0(X_\alpha)$ defined in \eqref{eqn:autom}.
The fixed points of~$\phi_1$ are~$\alpha^{-1}\pm\bigl(\sqrt{1-\alpha}\bigr)/\alpha$. Since $|\alpha|>1$, none of these roots is real; one lies above the real axis and the other below. The fundamental domain for~$\Gamma$ that we fixed in Lemma~\ref{thm:rep} lies at the left of the boundary geodesic going from~$\tau={\rm i}\infty$ to~$\tau=0$. This implies that the image, via~$t$, of the fundamental domain lies above the curve on~$X_\alpha$ that joins the punctures~$0$,~$\alpha^{-1}$. This implies that the root we have to choose is the one with positive imaginary part. If~$\tau_{1}=\hat{c}_{\phi_1}+1/\sqrt{D_1}$ is the fixed point on~$\Po$ of the lift of~$\phi_1$, we have
\[
t(\tau_{1}) = z_1 = \frac{1}{\alpha}+\frac{\sqrt{1-\alpha}}{\alpha}.
\]

Similar considerations apply to the choice of the fixed points~$z_2$ of the third non-trivial involution of~$X_\alpha$.
\end{proof}

\section{Finding the Fuchsian value}\label{sec:find}
\subsection{``Potential'' modular forms}\label{sec:cons}
The family of differential equations associated to~$X_\alpha$, determined from the general one in~\eqref{eqn:unifschw} by using the relations~\eqref{eqn:rel}, is given by
\[
\frac{{\rm d}^2y(t)}{{\rm d}t^2} + \left(\frac{P(t)'^2-P(t)P(t)''}{4P(t)^2}-\frac{t-m_0\alpha^{-1}}{2P(t)}\right)y(t)=0 ,\qquad P(t):=t(t-1)\big(t-\alpha^{1}\big) .
\]
In the following, we will not consider the above differential equation, but the projectively equivalent one
\begin{equation}\label{eqn:Heun}
\frac{{\rm d}^2y(t)}{{\rm d}t^2}+\frac{P'(t)}{P(t)}\frac{{\rm d}y(t)}{{\rm d}t}+\frac{(t-\rho)}{P(t)}y(t)=0,
\end{equation}
which is known as the~\emph{Heun equation}. The new accessory parameter~$\rho$ is related to~$m_0$ by~$m_0=1+\alpha-2\rho$.
As it will be clear, our results do not depend on the choice of the differential equation. We work with the Heun equation because its solutions are better behaved from the modular point of view; this will be relevant in numerical applications of our main result. It can be shown in fact that when~$\rho=\rho_{\rm F}$ is the Fuchsian value, the holomorphic solution lifts to (the square root of) a weight two modular form with a double zero in the cusp where the Hauptmodul~$t$ has its unique pole (see Appendix~\ref{appendixA} for the details).

The differential equation~\eqref{eqn:Heun} has at every finite singularity a~holomorphic solution and a~solution with a~logarithmic singularity. In particular, near the regular singular point~$t=0$ a~basis of solutions is given by
\begin{gather*}
y_{\alpha,\rho}(t)= \sum_{n\ge 0}{a_n(\alpha,\rho)t^n}= 1+\alpha\rho{t}+\frac{\alpha^2}{4}\big(\rho^2-2\rho(\alpha+1)-\alpha\big)t^2+\cdots, \\
\hat{y}_{\alpha,\rho}(t)= \log(t)y_{\alpha,\rho}(t) + \sum_{n\ge 0}{b_n(\alpha,\rho)t^n}= \log(t)y_{\alpha,\rho}(t)+\alpha(-2\rho+\alpha+1)t+\cdots ,
\end{gather*}
where the coefficients $a_n(\alpha,\rho)$, $b_n(\alpha,\rho)$ are polynomials in $\rho$ of degree~$n$ and satisfy the following linear recursions (Frobenius method)
\begin{gather*}
\alpha n^2a_{n-1}(\rho)-\big((\alpha+1)\big(n^2+n\big)+\rho\big)a_n(\rho)+(n+1)^2a_{n+1}(\rho)=0,
\\
\alpha n^2b_{n-1}(\rho)-\big((\alpha+1)\big(n^2+n\big)+\rho\big)b_n(\rho)+(n+1)^2b_{n+1}(\rho)\\
\qquad{} +2\alpha na_{n-1}(\rho)-(2n+1)(\alpha+1)a_n(\rho)+2(n+1)a_{n+1}(\rho)=0,
\end{gather*}
with initial data~$a_n=0$ if $n<0$, $a_0=1$ and~$b_n=0$ if~$n\le 0$.

The relevant function for the uniformization of $X_\alpha$ is the ratio of the two solutions~$y_\rho$,~$\hat{y}_\rho$. However, due to the logarithmic term, using power series it is more appropriate to work with the exponential of this ratio
\begin{equation}\label{eqn:Qd}
Q_{\alpha,\rho}(t):=\exp(\hat{y}_{\alpha,\rho}(t)/y_{\alpha,\rho}(t))=\sum_{n\ge 1}{Q_n(\alpha,\rho)t^n}= t+\alpha(-2\rho+\alpha+1)t^2+\cdots.
\end{equation}
The function $Q_{\alpha,\rho}(t)$ is a local biholomorphism as a function of $t$; inverting the series \eqref{eqn:Qd} we find the $Q$-expansion of its local inverse $t_{\alpha,\rho}(Q)$ around $Q=0$:
\begin{equation}\label{eqn:tQ}
t_{\alpha,\rho}(Q)=\sum_{n\ge1}{{t}_n(\alpha,\rho)Q^n}=Q-\alpha(-2\rho+\alpha+1)Q^2+\cdots.
\end{equation}
Finally, substituting the above series $t_{\alpha,\rho}(Q)$ into the holomorphic solution $y_{\alpha,\rho}(t)$, we get a new power series in $Q$:
\begin{gather}
f_{\alpha,\rho}(Q):=y_{\alpha,\rho}(t_{\alpha,\rho}(Q))=\sum_{n\ge 0}{{f}_n(\alpha,\rho)Q^n}\nonumber\\
\hphantom{f_{\alpha,\rho}(Q):=y_{\alpha,\rho}(t_{\alpha,\rho}(Q))}{}
= 1+\alpha\rho{Q}+\frac{\alpha^2}{4}\big(9\rho^2-2\rho(\alpha+1)-\alpha\big)Q^2+\cdots.\label{eqn:fQ}
\end{gather}
When the accessory parameter specializes to the Fuchsian value $\rho_{\rm F}$ the ratio $\hat{y}_{\alpha,\rho_{\rm F}}(t)/y_{\alpha,\rho_{\rm F}}(t)$ gives a coordinate on the universal covering~$\Po$ of~$X_\alpha$.
It follows from~\eqref{eqn:Qd} that~$Q_{\alpha,\rho_{\rm F}}(t)$ is a local parameter at the cusp~$\infty$ and that~$t_{\alpha,\rho_{\rm F}}(Q)$ is the local expansion of the Hauptmodul~$t\colon\Po\to X_\alpha$ in the parameter~$Q$.
A comparison between the expressions~\eqref{eqn:tQ} and~\eqref{eqn:expt} gives
\begin{equation}\label{eqn:newdefr}
Q=rq,\qquad\text{where}\quad q={\rm e}^{2\pi {\rm i}\tau}, \quad \tau\in\Po,
\end{equation}
for some non-zero $r \in \C$.
It follows that the $Q$-expansions \eqref{eqn:tQ}, \eqref{eqn:fQ} of $t_{\rho_{\rm F}}(Q)$ and $f_{\rho_{\rm F}}(Q)$ can be turned into $q$-expansions, which eventually make them holomorphic functions on $\Po$:
\begin{alignat*}{3}
& t(\tau):=t_{\alpha,\rho_{\rm F}}(rq)=\sum_{n\ge 1}{\hat{t}_nq^n},\qquad&& \hat{t}_n={t}_n(\alpha,\rho_{\rm F})r^n,&\\
& f(\tau):=f_{\alpha,\rho_{\rm F}}(rq)=\sum_{n\ge 0}{\hat{f}_nq^n},\qquad && \hat{f}_n={f}_n(\alpha,\rho_{\rm F})r^n.&
\end{alignat*}
From the discussion following~\eqref{eqn:Heun} we conclude that~$f^2$ is a weight two modular form with respect to the uniformizing group of~$X_\alpha$. On the contrary, the expansions~$t_{\alpha,\rho}(Q)$ and~$f_{\alpha,\rho}(Q)$ are \emph{``potential'' modular forms} in the sense that they extends to holomorphic functions on~$\Po$ with modular properties only for the correct value~$\rho_{\rm F}$ of~$\rho$. In the following we see them as functions depending on the parameter~$\rho$.

\subsection{``Potential'' cusp representatives}
Consider a four-punctured sphere~$X_\alpha$ where~$\alpha^{-1}$ is as in~\eqref{eqn:norma} and let~$z_j$, $j=0,1,2$ be the fixed points of the automorphisms of~$X_\alpha$ specified in Lemma~\ref{thm:fixed}.
We are going to consider a subset~$\mathcal{P}$ of the set of accessory parameters with a special property.
\begin{Definition}
For every~$\rho$ consider the power series~$t_{\alpha,\rho}(Q)$ defined in~\eqref{eqn:tQ} and let~$D_\rho$ denote its disk of convergence centered in~$Q=0$. We say that~$\rho\in\mathcal{P}$ if, for every~$j=0,1,2$, there exist~$Q_j\in D_\rho$ such that~$t_{\alpha,\rho}(Q_j)=z_j$.
\end{Definition}
This condition is not satisfied by most accessory parameters~$\rho$, but it is certainly satisfied by the Fuchsian parameter~$\rho_{\rm F}$ and, consequently, by an open subset of the set~$\mathcal{B}$ of parameters realizing a quasifuchsian uniformization of~$X_\alpha$.
For~$\rho\in\mathcal{P}$ and for~$j=0,1,2$ the function
\begin{equation*}
\frac{1}{t_{\alpha,\rho}(Q)-z_j}=\sum_{n\ge0}{{T}_{j,n}(\rho)Q^n},\qquad j=0,1,2,
\end{equation*}
has a simple pole in~$Q_j$ and is holomorphic in a punctured domain containing~$Q_j$. It follows that the limits
\[
 Q_j(\rho):= \lim_{n\to\infty}{\frac{{T}_{j,n}(\rho)}{{T}_{j,n+1}(\rho)}} ,\qquad j=0,1,2 ,
\]
exist for every fixed value of~$\rho\in\mathcal{P}$ and in fact define complex-valued functions of~$\rho$.
Finally, define the following real-valued functions of~$\rho\in\mathcal{P}$
\begin{equation*}
c_j(\rho):=\operatorname{Re}\left(\frac{1}{2\pi {\rm i}}\log\left(\frac{Q_j(\rho)}{Q_0(\rho)}\right)\right),\qquad j=0,1,2,
\end{equation*}
where $\log(z):=\log|z|+{\rm i}\operatorname{Arg}(z)$ and $\operatorname{Arg}(z)\in(-\pi,\pi]$.
The basic properties of the functions~$c_j(\rho)$ are given in the following lemma.
\begin{Lemma}\quad
\begin{enumerate}\itemsep=0pt
\item[$1.$] $c_0(\rho)= 0$, and $c_j(\rho)\in (-1/2,1/2]$ for every~$\rho\in\mathcal{P}$.
\item[$2.$] For every $\rho\in\mathcal{P}$ $c_i(\rho)\neq c_j(\rho)$ if~$i\neq j$.
\end{enumerate}
\end{Lemma}

In the next fundamental proposition we attach a Fuchsian group to every differential equation~\eqref{eqn:Heun} with~$\rho\in\mathcal{P}$. We remark that the Fuchsian group~$\Gamma(\rho)$ attached to~$\rho\in\mathcal{P}$ is in general not the monodromy group of the associated differential equation, which is a Kleinian non Fuchsian group, but a group constructed by considering the automorphisms of the four-punctured sphere. The group~$\Gamma(\rho)$ is the monodromy group only when~$\rho=\rho_{\rm F}$ is the Fuchsian parameter (point~3 of the proposition).
\begin{Proposition}\label{thm:fund}\quad\samepage
\begin{enumerate}\itemsep=0pt
\item[$1.$] For every $\rho\in\mathcal{P}$ there exist a unique torsion-free Fuchsian group~$\Gamma(\rho)$ of genus zero with four cusps and nonequivalent cusp representatives $0$, $c_1(\rho)$, $c_2(\rho)$, $\infty$.
\item[$2.$] For every fixed $\rho\in\mathcal{P}$ let~$x_0$, $x_1$, $x_2$ be real numbers such that
\[\{x_0,x_1,x_2\}=\{0,c_1(\rho),c_2(\rho)\}\qquad\text{and}\qquad x_0<x_1<x_2.
\]
Define, for $j=0,1,2$,
\[
S_j=S_j(\rho):=\begin{pmatrix}
1+x_jA_j & -x_j^2A_j\\
A_j & 1-x_jA_j
\end{pmatrix},
\]
where~$A_j=A_j(\rho):= (x_j-x_{j-1} )^{-1} (x_{j+1}-x_j )^{-1}$, and~$x_{-1}:=x_2-1$, $x_3:=x_0+1$. Then
\[
\Gamma(\rho)=\langle T,S_0(\rho),S_1(\rho),S_2(\rho)\rangle \qquad\text{and}\qquad S_2S_1S_0T^{-1}=1 .
\]
\item[$3.$] When $\rho=\rho_{\rm F}$ is the Fuchsian value, the group $\Gamma(\rho_{\rm F})$ is the uniformizing group of~$X_\alpha$.
\end{enumerate}
\end{Proposition}

\begin{proof}Consider three real numbers~$x_0,x_1,x_2\in(-1/2,1/2]$ such that~$x_0<x_1<x_2$. We shall associate to the triple~$(x_0,x_1,x_2)$ a torsion-free Fuchsian group with four cusps whose representatives are~$x_0$, $x_1$, $x_2$ and $\infty$ by using Poincar\'e's theorem.

Let $x_1':=\frac{x_2+x_1(x_2-1)-x_0x_1}{x_2-x_0}$. Using the properties~$-1/2<x_j\le1/2$ and~$x_0<x_1<x_2$, it is easy to verify that $x_2<x_1'<x_0+1$.
Consider~$\Po$ as a model of the hyperbolic plane, and let~$\mathcal{F}\subset\Po$ be the hyperbolic geodesic polygon with vertices~$\{x_0,x_1,x_2,x_1',x_0+1,\infty\}$. A simple calculation shows that the set of transformations $G:=\big\{T,S_2,S_2S_1^{-1}\big\}$ is a side-pairing for the geodesic boundary of $\mathcal{F}$ and~$S_2S_1S_0T^{-1}=1$. We can conclude by Poincar\'e's theorem~(see~\cite{Beardon}) that the group generated by the transformations in~$G$ is a Fuchsian group of genus zero with no torsion and four cusps and with fundamental domain~$\mathcal{F}$.
The first two points of the proposition follow by choosing $x_j=x_j(\rho)$ as in the statement. We denote by $\Gamma(\rho)$ the Fuchsian group obtained in this way.

We prove point 3. When~$\rho=\rho_{\rm F}$ we know by~\eqref{eqn:newdefr} that~$Q=r{\rm e}^{2\pi {\rm i}\tau}$ for some non-zero~$r\in\C$.
It follows that~$Q_j=r{\rm e}^{2\pi {\rm i}\tau_j}$, $j=0,1,2$, where~$\tau_j$ is the fixed point in~$\Po$ of the lifting of the automorphism~$\phi_j$ of~$X_\alpha$ (see Section~\ref{sec:cusprep}). Using the description of~$\tau_j$ in~\eqref{eqn:fph} we see that
\[
\log(Q_j/Q_0)/(2\pi {\rm i}) = \tau_j -\tau_0 = {c}_j +{\rm i}\big(1/\sqrt{D_j}-1/\sqrt{D_0}\big),\qquad j=1,2,
\]
where~$c_1$, $c_2$ are inequivalent cusps of the uniformizing group of~$X_\alpha$. Since~$c_0(\rho_{\rm F})=0$ and $c_j(\rho_{\rm F})=c_j$, $j=1,2$, we conclude that the group~$\Gamma(\rho_{\rm F})$ constructed in point~2 is the uniformizing group of~$X_\alpha$.
\end{proof}

Similarly to the ``potential'' modular forms of Section~\ref{sec:cons}, the functions~$c_j(\rho)$ are cusps of the uniformizing group of~$X_\alpha$ for the value~$\rho_{\rm F}$ of the accessory parameter. For this reason, we call the~$c_j(\rho)$ \emph{``potential'' cusp representatives}, even though they are actually cusps for the group~$\Gamma(\rho)$ for every~$\rho\in\mathcal{P}$.

\subsection{Finding the Fuchsian value}\label{sec:findfin}
In the previous section we constructed a Fuchsian group~$\Gamma(\rho)$ from the differential equation~\eqref{eqn:Heun} attached to the four-punctured sphere~$X_\alpha$ if~$\rho\in\mathcal{P}$. When~$\rho=\rho_{\rm F}$ the group~$\Gamma(\rho_{\rm F})$ is the uniformizing group of~$X_\alpha$ and the function~$t_{\alpha,\rho_{\rm F}}(rq)$ obtained by inverting the exponential of the ratio of independent solutions of~\eqref{eqn:Heun} is a modular function with respect to~$\Gamma(\rho_{\rm F})$.

The idea is to mimic this situation in the case~$\rho\in\mathcal{P}$ in order to check whether~$\rho=\rho_{\rm F}$ by checking the modularity of~$t_{\alpha,\rho}(Q)$ with respect to~$\Gamma(\rho)$. To do this, we need to make~$t_{\alpha,\rho}(Q)$ a~function on~$\Po$ for every~$\rho\in\mathcal{P}$. We do it as follows.

In the proof of Proposition~\ref{thm:fund} we showed that when~$\rho=\rho_{\rm F}$ we have~$Q_0=r{\rm e}^{2\pi {\rm i}\tau_0}$ for some non-zero~$r\in\C$. Moreover~$\tau_0={\rm i}/\sqrt{D_0}={\rm i} (c_1(1-c_2) )^{1/2}$ as follows from~\eqref{eqn:fph} and~\eqref{eqn:Di}. We can then easily determine~$r=\log(Q_0)/\exp\big({-}2\pi/\sqrt{D_0}\big)$.
For a generic~$\rho\in\mathcal{P}$ then it makes sense to define
\begin{equation*}
r(\rho):=\frac{Q_0(\rho)}{\exp\big({-}2\pi/\sqrt{D_0(\rho)}\big)},\qquad D_0(\rho)=\frac{1}{c_1(\rho) (1-c_2(\rho) )},
\end{equation*}
and make~$t_{\alpha,\rho}(Q)$ into a holomorphic function on~$\Po$ by setting
\[
t_{\alpha,\rho}(\tau):=t_{\alpha,\rho}\bigl(r(\rho){\rm e}^{2\pi {\rm i}\tau}\bigr) ,\qquad \tau\in\Po .
\]
The function~$t_{\alpha,\rho}(\tau)$ is modular with respect to the group~$\Gamma(\rho)$ if the following functions
\begin{equation}
\label{eqn:modfin}
E_{\alpha,j}(\rho,\tau):=t_{\alpha,\rho}\bigl(r(\rho) {\rm e}^{2\pi {\rm i}S_j(\rho)\tau}\bigr)-t_{\alpha,\rho}\bigl(r(\rho){\rm e}^{2\pi {\rm i}\tau}\bigr),\qquad j=0,1,2 ,
\end{equation}
where~$S_j(\rho)$ are the generators of~$\Gamma(\rho)$, are zero for every~$\tau$ in a fundamental domain for~$\Gamma(\rho)$.

\begin{Theorem}\label{thm:main}Let $X_\alpha$ be a four-punctured sphere and let~$E_{\alpha,j}(\rho,\tau)$, $j=0,1,2$ be as in~\eqref{eqn:modfin}. The Fuchsian value $\rho_{\rm F}$ for the uniformization of $X_\alpha$ is the unique zero of the system of equations
\[
 E_{\alpha,j}(\rho,\tau)=0,\qquad j=0,1,2,
\]
for every $\tau$ in a fundamental domain of~$\Gamma(\rho)$.
\end{Theorem}
\begin{proof}It is clear that~$\rho_{\rm F}$ is a zero of~$E_{\alpha,j}(\rho,\tau)$ for~$j=0,1,2$ and every~$\tau\in\Po$.

Let~$\rho_1\in\mathcal{P}$ be such that the identity in the statement holds for every~$j=0,1,2$. Since the function~$t_{\alpha,\rho_1}(Q)$ is univalent in~$Q$ and~$t_{\alpha,\rho_1}(0)=0$ it follows that~$t_{\alpha,\rho_1}(\tau)$ is never zero on~$\Po$ and then, being holomorphic, it is a Hauptmodul for the group~$\Gamma(\rho_1)$. The only thing to check is that~$\Gamma(\rho_1)=\Gamma(\rho_{\rm F})$ and that~$t_{\alpha,\rho_1}(\tau)$ is a covering map for~$X_\alpha$.

Since~$t_{\alpha,\rho_1}(\tau)$ has a simple pole at one cusp it maps the Riemann surface~$\Po/\Gamma(\rho_1)$ to the punctured sphere~$\pro^1\smallsetminus\big\{\infty,0,a,b^{-1}\big\}$ for some~$a,b\in\C\smallsetminus\{0\}$, $a\neq b^{-1}$.
We can assume that~$t_{\alpha,\rho_1}(\tau)$ maps the cusp~$0$ to~$b^{-1}$. It follows that~$a^{-1}t_{\alpha,\rho_1}(\tau)$ is a Hauptmodul for the punctured sphere~$X_{ab}=\pro^1\smallsetminus\big\{\infty,0,1,a^{-1}b^{-1}\big\}$ and that~$a^{-1}t_{\alpha,\rho_1}(\tau)$ is normalized as in Lemma~\ref{thm:norm}. We can then obtain the expansion at~$\infty$ of~$a^{-1}t_{\alpha,\rho_1}(\tau)$ from a basis of solutions~$\{y_{ab,\hat{\rho}}(t),\hat{y}_{ab,\hat{\rho}}(t)\}$ of the uniformizing differential equation of~$X_{ab}$ as in Section~\ref{sec:cons}, i.e.,
\[
a^{-1}t_{\alpha,\rho_1}(\tau)=t_{ab,\hat{\rho}}(\tau)=\sum_{n\ge0}{t_n(ab,\hat{\rho})r(\hat{\rho})^nq^n} ,
\]
where~$\hat{\rho}$ is the Fuchsian parameter associated to the uniformizing differential equation of~$X_{ab}$.
On the other hand, we can describe~$t_{\alpha,\rho_1}(\tau)$ at~$\infty$ with a power series constructed from a basis~$\{y_{\alpha,\rho_1}(t),\hat{y}_{\alpha,\rho_1}(t)\}$ of solutions of the differential equation on~$X_\alpha$ with accessory parameter~$\rho_1$
\[
t_{\alpha,\rho_1}(\tau)=\sum_{n\ge0}{t_n(\alpha,\rho_1)r(\rho_1)^nq^n} .
\]
Finally, since~$\{y_{ab,\hat{\rho}}(t),\hat{y}_{ab,\hat{\rho}}(t)\}$ is a basis of solutions of the uniformizing equation for~$X_{ab}$, and by comparing the power series representations of~$t_{\alpha,\rho_1}$ we get
\begin{align*}
a\exp\left(\frac{\hat{y}_{ab,\hat{\rho}}}{y_{ab,\hat{\rho}}}\right)&=\left(a\sum_{n=1}^\infty{t_n(ab,\hat{\rho})r(\hat{\rho})^nq^n}\right)^{-1} =\left(\sum_{n=1}^\infty{{t_n(\alpha,\rho_1)}\frac{r(\rho_1)^n}{r(\hat{\rho})^n}r(\hat{\rho})^nq^n}\right)^{-1}\\ &=\frac{r(\rho_1)}{r(\hat{\rho})}\exp\left(\frac{\hat{y}_{\alpha,\rho_1}}{y_{\alpha,\rho_1}}\right) ,
\end{align*}
where~$^{-1}$ denotes the compositional inverse. It follows that the ratio of solutions of the differential equation for~$X_{\alpha}$ with accessory parameter~$\rho_1$ and of the uniformizing one of~$X_{ab}$ differ only by a constant factor. This implies that the ratio~$\hat{y}_{\alpha,\rho_1}/y_{\alpha,\rho_1}$ induces a biholomorphism~$\widetilde{X}_{\alpha}\rightarrow\Po$, i.e., that~$\rho_1$ is the Fuchsian parameter. By the uniqueness of the Fuchsian parameter we can conclude that~$\rho_1=\rho_{\rm F}$.
\end{proof}

\section{Example: local expansion of the Fuchsian value function}\label{sec:xpl}
\subsection{Numerical computation of the Fuchsian value}\label{sec:method}
We first explain how to use Theorem~\ref{thm:main} to approximate numerically the Fuchsian value for a~given four-punctured sphere~$X_\alpha=\pro^1\smallsetminus\big\{\infty,1,0,\alpha^{-1}\big\}$.
The behavior of the Fuchsian value under the action of the anharmonic group (the group of order six generated by~$z\to z^{-1}$ and~$z\to 1-z$) and the action of complex conjugation is known~\cite{KRV}; we need then to consider only the case\looseness=-1
\[
\alpha^{-1}\in \{z\in\C\colon 0\le\operatorname{Re}(z)\le1/2,\,|z|\le1\}\smallsetminus\{0\} .
\]
It follows from Theorem~\ref{thm:main} that to compute the Fuchsian value for~$X_{\alpha}$ is equivalent to compute the common zero of the equations~$E_{\alpha,j}(\rho,\tau)=0$ in~\eqref{eqn:modfin}.
Notice that all the quantities involved in the definition of~$E_{\alpha,j}(\rho,\tau)$ are computable, as functions of~$\rho$, from the Frobenius solutions of the differential equation~\eqref{eqn:Heun}.
We proceed as follows. Fix~$\tau_0\in\Po$ and consider, for~$j=0$, the equation~$E_{\alpha,0}(\rho,\tau_0)=0$. We use Newton's method to find the zero of this equation that is the Fuchsian value. As it is well known, the Newton method works if we are able to give an initial guess for the zero that it is close enough to it. In other words, to start the iteration we should choose a value of~$\rho$ that is close to the Fuchsian value.
Since the function associating to a~four-punctured sphere its Fuchsian value is continuous, a good choice for the initial value is the Fuchsian value of a four-punctured sphere~$X_\beta$ with~$\beta^{-1}$ close to~$\alpha^{-1}$.
There are four exceptional choices of~$\beta^{-1}$ for which is possible to determine exactly the value of the Fuchsian parameter via symmetries (see \cite[Section~7]{Hempel} and~\cite{ZagApery}; the uniformizing groups in these cases are conjugated to congruence subgroups of~${\rm SL}_2(\Z)$). These choices and their Fuchsian value are displayed in the following table
$$
\begin{array}{lcccc}
\toprule
\beta^{-1} & \frac{1}{2} & \frac{1}{9} & \frac{25-11\sqrt{5}}{50} & \frac{1+{\rm i}\sqrt{3}}{2} \\
\midrule
\rho_{\rm F} & 1 & 3 & \frac{35+15\sqrt{5}}{2} & \frac{3-{\rm i}\sqrt{3}}{6} \\
\bottomrule
\end{array}
$$
Assume for now that the value of~$\alpha^{-1}$ is close enough to one of the special values~${\beta}^{-1}$. In this case we start the Newton's method with the Fuchsian value~$\rho_{\rm F}\big({\beta}^{-1}\big)$. The iteration gives the approximation of a zero of~$E_{\alpha,0}(\rho,\tau_0)$ that is close to the Fuchsian parameter. To verify if it really is the Fuchsian parameter, we check whether it is a zero of other equations~$E_{\alpha,j}(\rho,\tau)=0$ for different choices of~$j$ and~$\tau\in\Po$.

In the case~$\alpha^{-1}$ is not close enough to one of the special values, one can gradually approach the computation of the Fuchsian value of~$X_\alpha$ by computing the Fuchsian value of some points between one of the special values~$\beta^{-1}$ and~$\alpha^{-1}$.

\subsection{Application: local expansion of the Fuchsian value function}
As an application, we compute numerically the local expansion of the function that associates to a four-punctured sphere $X_\alpha$ the Fuchsian value $\rho_{\rm F}(X_\alpha)$. We can define this function in full generality for an~$n$-punctured sphere~$X$. Let~$W_n=\{(w_1,\dots, w_{n-3})\,|\,w_i\neq w_j\text{ if }i\neq j,\,w_i\neq 0,1\}$ and consider the function
\[
\rho_{\rm F}\colon \ W_{n}\to\C^{n-3},\qquad w=(w_1,\dots,w_{n-3})\mapsto \rho_{\rm F}(w)=(\rho_1,\dots,\rho_{n-3}),
\]
where $\rho_{\rm F}=(\rho_1,\dots, \rho_{n-3})$ is the Fuchsian value for the $n$-punctured sphere
\[
X=\pro^1\smallsetminus\{w_1,w_2,\dots,w_{n-3},0,1,\infty\}.
\]
Kra \cite{Kra} proved that the function~$\rho_{\rm F}$ is real-analytic, but non complex-analytic; in particular, if~$z$ is a local parameter on~$W_n$, the function $\rho_{\rm F}$ has a local expansion near the point $z_0\in W_n$ of the form
\begin{equation}\label{eqn:exp}
\rho_{\rm F}(z_0+z)=\sum_{j,k\ge 0}{a_{j,k}z^j\bar{z}^k},\qquad a_{j,k}\in\C.
\end{equation}
In the following we concentrate on the case $n=4$ and study the expansion of the function
\[
\rho_{\rm F}\colon \ \C\setminus\{0,1\}\to\C,\qquad \alpha^{-1}\mapsto\rho_{\rm F}(X_\alpha),
\]
around the points $\alpha^{-1}=1/2,1/3$, and~$2/5+3{\rm i}/10$.

The computation of the coefficients of~\eqref{eqn:exp} goes as follows.
Fix~$\alpha^{-1}\in\C\smallsetminus\{0,1\}$ and consider for every~$m\in\N$ the line~$L_m$: $\operatorname{Im}(z)=\operatorname{Re}(z)/m-\alpha^{-1}/m$. The expansion of~$\rho_{\rm F}$ on each line~$L_m$ depends only on one real variable~$x$, since~$z\in L_m$ if~$z=\alpha^{-1}+x(1+{\rm i}/m)$ and then
\begin{gather}\label{eqn:alfa}
\rho_{\rm F}\big(\alpha^{-1}+z\big)=\sum_{j,k\ge 0}{a_{j,k}}(1+{\rm i}/m)^j(1-{\rm i}/m)^kx^{k+j}=\sum_{s=0}^\infty{b_s(m)x^m} ,\qquad b_s(m)\in\C .
\end{gather}
The coefficients~$b_s(m)$ are easily computed once we know enough values of the function~$\rho_{\rm F}$ on the line~$L_m$ near~$\alpha^{-1}$. These can be computed by using the method illustrated in Section~\ref{sec:method} and the computer algebra system~PARI.
The coefficients~$a_{j,k}$ in the expansion of~$\eqref{eqn:exp}$ are then obtained by computing the expansion of~$\rho_{\rm F}$ along the lines $L_m$ for different values of~$m\in \N$ (the number of those depends on the number of $a_{j,k}$ one wants to compute) and exploiting the relation between~$a_{j,k}$ and~$b_s(m)$ in~\eqref{eqn:alfa}. The result of the computations for~$\alpha^{-1}=1/2,1/3,2/5+3{\rm i}/10$ are given in Tables~\ref{table1},~\ref{table2}, and~\ref{table3} respectively (the numbers in the tables are approximations of the actual values).
\begin{table}[h!]\centering
\caption[]{$\alpha^{-1}=1/2$, $\rho_{\rm F}=1$.}\label{table1}
$\begin{array}{lr}
\toprule
a_{1,0}=-1.2311296972& a_{0,1}=0.0638754899 \\
\midrule
a_{2,0}=2.4622593944& a_{1,1}=-0.1277509798\\
a_{0,2}=0 & \\
\midrule
a_{3,0}=-4.8236918585& a_{2,1}=0.1890620793\\
a_{1,2}=0.0117490877& a_{0,3}=0.0630206931\\
\midrule
a_{4,0}=9.6473837171&a_{3,1}=-0.3781241587\\
a_{2,2}=-0.0234981755&a_{1,3}=-0.1260413862\\
a_{0,4}=0 &\\
\midrule
a_{5,0}=-19.094665845& a_{4,1}=0.6673769276\\
a_{3,2}=0.0466379026& a_{2,3}=0.1888625099\\
a_{1,4}=0.0233189513& a_{0,5}=0.1334753855\\
\bottomrule
\end{array}$\vspace{-4mm}
\end{table}
\begin{table}[ht!]\centering
\caption[]{$\alpha^{-1}=1/3$, $\rho_{\rm F}=1.29101$.}\label{table2}
$\begin{array}{lr}
\toprule
a_{1,0}=-2.711485382 & a_{0,1}=0.1025201219 \\
\midrule
a_{2,0}=8.0641055547 & a_{1,1}=-0.2750330946 \\
a_{0,2}=-0.0606264558 & \\
\midrule
a_{3,0}=-23.9531822161 & a_{2,1}= 0.6879078089 \\
a_{1,2}=0.1854950416& a_{0,3}=0.1686761471\\
\midrule
a_{4,0}=71.4914941489&a_{3,1}=-1.9281043695 \\
a_{2,2}=-0.4821058340&a_{1,3}=0.4798627522\\
a_{0,4}=0.2922654268&\\
\midrule
a_{5,0}=-213.5180837271 & a_{4,1}=5.4375565883\\
a_{3,2}=1.3699401215& a_{2,3}=1.2274288636\\
a_{1,4}=0.8681638848& a_{0,5}=0.7367927177\\
\bottomrule
\end{array}$\vspace{-4mm}
\end{table}
\begin{table}[ht!]\centering
\caption[]{$\alpha^{-1}=\frac{4+3{\rm i}}{10}$, $\rho_{\rm F}=0.86175-0.38528 {\rm i}$}\label{table3}
$\begin{array}{lr}
\toprule
a_{1,0}=-0.3328603817+1.2004121803{\rm i} & a_{0,1}=0.0512782931-0.0256391465{\rm i} \\
\midrule
a_{2,0}=-0.8635602303-2.2690958807{\rm i}& a_{1,1}=-0.0638839260+0.0719692197{\rm i}\\
a_{0,2}=-0.0223137191-0.03022272424{\rm i}\\
\midrule
a_{3,0}=4.0817033346+2.5852542694{\rm i} & a_{2,1}=-0.0024129061-0.1684860963{\rm i}\\
a_{1,2}=0.0599341643+0.0300072676{\rm i}& a_{0,3}=0.01118030554+0.0173932256{\rm i}\\
\midrule
a_{4,0}=-9.5604869979+0.7791704955{\rm i}&a_{3,1}=0.1800042196+0.2360936333{\rm i}\\
a_{2,2}=-0.1150333973+0.0337939365{\rm i}&a_{1,3}=-0.0337943677-0.0275496197{\rm i}\\
a_{0,4}=-0.0003411753-0.0596118160{\rm i}\\
\midrule
a_{5,0}=14.2865112529-12.6279958566{\rm i}& a_{4,1}=-0.5419356177-0.1625218437{\rm i}\\
a_{3,2}=0.1261663672-0.1720412035{\rm i}& a_{2,3}=0.0798727692-0.0023875661{\rm i}\\
a_{1,4}=0.0384447463+0.0991919592{\rm i}& a_{0,5}=-0.0711856680+0.0428083431{\rm i}\\
\bottomrule
\end{array}$\vspace{-4mm}
\end{table}

A first interesting observation we can make from the numerical data is about the size of the coefficients~$a_{j,k}$. For every~$s=1,\dots,5$ we have that~$|a_{s,0}|>|a_{j,k}|$ with~$j+k=s$. In other words, the holomorphic part of the expansion of~$\rho\big(\alpha^{-1}\big)$ seems to be larger than the rest. This suggests that the Fuchsian parameter function may be quasiregular, i.e., it may satisfy the inequality
\[
\frac{\partial{\rho_{\rm F}}}{\partial\bar{z}} \le k\frac{\partial\rho_{\rm F}}{\partial z}
\]
for some~$k<1$. This would in particular imply that the Fuchsian parameter map is open.

We further notice the following relations among the coefficients in Table~\ref{table1}:
\begin{gather}\label{eqn:boh1}
a_{2,0}=-2a_{1,0},\qquad a_{1,1}=-2a_{0,1},\\
a_{0,2}=0=a_{0,4},\qquad
a_{2,1}=3a_{0,3},\\
\label{eqn:boh4}
a_{3,2}=2a_{1,4},\qquad a_{4,1}=5a_{0,5}.
\end{gather}
These numerical identities can be proven by using the symmetry of $\rho_{\rm F}$ near $1/2$, and a result of Takhtajan and Zograf \cite{TZ1}. Analogous identities in Table~\ref{table2} or Table~\ref{table3} can be proved similarly.
The point $z_0=1/2$ is the fixed point of the involution $z\mapsto 1-z$. It is known (see~\cite{KRV}) that the following identity holds
\[
\rho_{\rm F}(1-z) = \frac{z\rho_{\rm F}(z)-1}{z-1}.
\]
It follows that, near the point $z_0=1/2$, one has
\[
(z-1/2)\rho_{\rm F}(1/2-z)=(1/2+z)\rho_{\rm F}(1/2+z)-1,
\]
which gives
\begin{equation}\label{eqn:uno}
\sum_{j,k\ge 0}{a_{j,k}\big[1+(-1)^{j+k}\big]z^j\overline{z}^k}+2\sum_{j,k\ge 0}{a_{j,k}\big[1-(-1)^{j+k}\big]z^{j+1}\overline{z}^k}-2=0.
\end{equation}
The above relation implies that
\begin{gather*}
\begin{split}
& a_{0,2k}=0\qquad\text{if } k\ge 1, \\
& a_{j+1,k}=-2a_{j,k}\qquad\text{if }k+j\text{ is odd}.
\end{split}
\end{gather*}
This explains why $a_{0,2}=a_{0,4}=0$.
The result of Takhtajan and Zograf \cite[formula~(4.1)]{TZ1} reduces to the following identity in the case of four-punctured spheres\footnote{this equation is formulated in~\cite{TZ1} in terms of the accessory parameters~$m_i$ appearing in the Schwarzian differential equation~\eqref{eqn:unifschw}; here we express it in terms of the accessory parameter of the Heun equation~\eqref{eqn:Heun}.}
\begin{equation}\label{eqn:de3}
(1-2z)\overline{\left(\frac{\partial\rho_{\rm F}}{\partial\overline{z}}\right)} = (1-2\overline{z})\frac{\partial\rho_{\rm F}}{\partial\overline{z}}.
\end{equation}
The differential equation~\eqref{eqn:de3} implies the following relations between the coefficients of the local expansion of~$\rho_{\rm F}$:
\[
(k+1)a_{j,k+1}-2ka_{j,k}=(j+1)a_{k,j+1}-2ja_{k,j},\qquad j,k\ge 0.
\]
It is easy to check that the relations \eqref{eqn:boh1}--\eqref{eqn:boh4} come from this one and from \eqref{eqn:uno}. For instance, by choosing $(j,k)=(0,1)$ in the identity above we get
\[
2a_{0,2}-2a_{0,1}=a_{1,1}.
\]
This, together with $a_{0,2}=0$, gives the first identity in \eqref{eqn:boh1}.

\appendix

\section{Modular derivation of the uniformizing differential equation}\label{appendixA}
Denote by~$\Gamma\subset{\rm SL}_2(\R)$ a genus zero Fuchsian group with no torsion and with~$n\ge 3$ inequivalent cusps. Normalize it by assuming that one of its cusps is at~$\infty$, and that this cusp has width one. Let~$t$ be a Hauptmodul and, without loss of generality, assume that its unique pole is at a~cusp~$c_0$ and its unique zero is at~$\infty$. Finally, let~$X_\Gamma$ be the~$n$-punctured sphere isomorphic to~$\Po/\Gamma$ via~$t$
\[
t\colon \ \Po/\Gamma \overset{\sim}{\rightarrow} X_\Gamma=\pro^1 \smallsetminus \{\alpha_1,\alpha_2,\dots,\alpha_{n-1}=0 ,\alpha_n=\infty\},
\]
where $\alpha_i\in\C\smallsetminus\{0,1\}$, $i=1,\dots,n-2$, and~$\alpha_i\neq\alpha_j$ if~$i\neq j$.

In the following we compute the differential equation satisfied by a certain modular form~$f$ with respect to~$t$. This differential equation is projectively equivalent to the differential equation~\eqref{eqn:unifschw} associated to the uniformization of~$X_\Gamma$, and in the case~$n=4$ reduces to the Heun equation~\eqref{eqn:Heun} considered in Section~\ref{sec:find}. In particular, this gives a purely modular definition of the accessory parameters, as it will be clear from the proof of Proposition~\ref{thm:f}.

Since the differential equation in~\eqref{eqn:unifschw} has order two it would be natural, according to Section~\ref{sec:funcor}, to consider a weight one modular form on $\Gamma$, which satisfies a second order differential equation. It is known however that not every group~$\Gamma$ admits weight one modular forms (for now we are only assuming that~$\Gamma$ is of genus zero and torsion free). It makes sense then to consider a square root of a modular form of weight two, since $\dim{M_2(\Gamma)}=n-1$ for every torsion-free genus zero group~$\Gamma$ with~$n$ cusps.
We choose to work with a weight two modular form whose zeros are concentrated in a certain cusp; as the next lemma shows, this choice is always possible.

\begin{Lemma}\label{thm:choice}
Let $\Gamma$ be torsion free and of genus zero, let $t$ be an Hauptmodul, and denote by~$c_0$ the cusp of~$\Gamma$ where~$t$ has its unique pole.
There exists a modular form $f\in M_2(\Gamma)$, unique up to scalar multiplication, with all its zeros in~$c_0$. In particular, $f$ has no zeros in~$\Po$.
\end{Lemma}
\begin{proof}Let $g\in M_2(\Gamma)$ and let $\sigma\in{\rm SL}(2,\R)$ be such that $\sigma{c_0}=\infty$. Let
\[
\bigl(g\bigr|_2\sigma^{-1}\bigr)(\tau)=\sum_{m\ge 0}{g_mq^m}
\]
denote the Fourier expansion of $g$ at $c_0$, where $q={\rm e}^{2\pi {\rm i}\tau/h},\tau\in\Po$, is a local parameter. It is known that the degree of the divisor associated to any $g\in M_2(\Gamma)$ is $d=n-2$.
Let $\phi$ be the map
\[
\phi\colon \ M_2(\Gamma)\to\C^d,\qquad g\mapsto (g_0,g_1,\dots,g_{d-1}).
\]
that sends a modular form of weight $2$ to the vector defined by its first $d$ Fourier coefficients at the cusp $c_0$. This map is linear.

The dimension of $M_2(\Gamma)$ is $n-1=d+1$, so the map $\phi$ has a non-trivial kernel of dimension~$\ge1$.
Let $f\in\operatorname{Ker}(\phi)$. Such $f$ can have at most $d$ zeros in ${\Po}\cup\{\text{cusps}\}$, and they are all in $c_0$ by construction.
Finally, let $f,g\in\operatorname{Ker}(\phi)$ be linearly independent. The ratio $f/g$ is a weight zero modular form holomorphic in $\Po$ and in all the cusps, since $f$ and~$g$ have all their zeros at the same cusp~$c_0$. This implies that~$f/g$ is a constant, i.e.,~$\dim\operatorname{Ker}(\phi)=1$.
\end{proof}

Given~$f$ and~$t$ as in Lemma \ref{thm:choice} we can construct all the modular forms of even weight on~$\Gamma$.
\begin{Lemma}\label{thm:dime}
Let $k\ge0$ be an integer, and let $f$ and $t$ be as in Lemma~{\rm \ref{thm:choice}}. The functions
\[
f^kt^i,\qquad i=0,\dots,k(n-2),
\]
form a basis of the space $M_{2k}(\Gamma)$.
\end{Lemma}
\begin{proof}By construction, the weight $2k$ modular form $f^k$ has $k(n-2)$ zeros at the cusp $c_0$ where $t$ has a simple pole, and these are the only zeros of~$f^k$.
It follows that $f^kt^i$ is a holomorphic modular form for every $i=0,\dots,k(n-2)$, and meromorphic for every other value of~$i$.
By looking at the location of the zeros, we can prove that the holomorphic functions in the statement are linearly independent. From the dimension formula for $M_{2k}(\Gamma)$ (see for example~\cite[Chapter~2]{Miyake}) we conclude that they form a basis.
\end{proof}

The second order linear differential operator associated to a square root of the modular form~$f\in M_2(\Gamma)$ in Lemma~\ref{thm:choice} and to the Hauptmodul~$t$ is given in the next proposition.

\begin{Proposition}\label{thm:f} Let $\Gamma$ be a genus zero torsion-free Fuchsian group with $n\ge 3$ cusps, and let~$t$ be a Hauptmodul such that $t\colon\Po/\Gamma\overset{\sim}{\rightarrow}X_\Gamma=\pro^1\smallsetminus\{\alpha_1,\dots,\alpha_{n-2},\alpha_{n-1}=0,\alpha_n=\infty\}$.
Denote by~$c_0$ the cusp of~$\Gamma$ where $t$ has its unique pole, and let $f\in M_2(\Gamma)$ be such that all its zeros are at the cusp $c_0$.
Then the differential operator~$L$ associated to a square root of~$f$ and to~$t$ is given by\looseness=-1
\begin{equation}\label{eqn:f}
L = \frac{{\rm d}}{{\rm d}t}\left(P(t)\frac{{\rm d}}{{\rm d}t}\right)+\sum_{i=0}^{n-3}\rho_it^i,
\end{equation}
where $P(t)=\prod_{j=1}^{n-1}{(t-\alpha_j)}$, $\rho_{n-3}=(n/2-1)^2$, and $\rho_0,\dots,\rho_{n-4}\in\C$ are uniquely determined by~$f$,~$t$.
\end{Proposition}

\begin{proof}Recall from Section~\ref{sec:funcor} that $L$ can be computed in terms of Rankin--Cohen brackets of~$f$ and~$t$ by
\begin{equation}\label{eqn:eqproof1}
L= \frac{{\rm d}^2}{{\rm d}t^2}+\frac{[f,t']_1}{2ft'^2}\frac{{\rm d}}{{\rm d}t}-\frac{[f,f]_2}{12f^2t'^2}.
\end{equation}
We have to write the coefficients of $L$ as rational functions of~$t$.
First we prove that
\[
(-1)^n\left(\prod_{j=1}^{n-2}{\alpha_j}\right)t'=fP(t).
\]
The ratio $t'/f$ is a meromorphic modular function, so it is a rational function of $t$. From the assumption on the zeros of $f$ it follows that the modular function $t'/f$ has a simple zero at every cusp different from~$c_0$, i.e., $n-1$ simple zeros (since these are the zeros of~$t'$). It has also a~unique pole of order $n-1$ at $c_0$, since~$f$ has $n-2$ zeros there and~$t'$ a~simple pole. The rational functions of $t$ with these zeros and poles are given by the polynomials~$\kappa^{-1}P(t)$, $\kappa\in \C^*$, where~$P(t)$ is as in the statement. Looking at the first coefficient of the $q$-expansion of~$t'/f$ at~$\infty$, we find the correct factor $\kappa=(-1)^{n}\prod_{j=1}^{n-2}{\alpha_j}$.

Next, we compute the brackets $[f,t']$, and $[f,f]_2$.
The first one is very easy
\begin{align*}
[f,t'] & =2ft''-2f't'=f(fP(t)\kappa^{-1})'-2f't' \\
& =2f't'+2\kappa^{-1} f^2P'(t)t'-2f't'=2\kappa^{-1} f^2t'P'(t).
\end{align*}
Dividing~$[f,t']$ by $2ft'^2=2\kappa f^2P(t)t'$ we see that the coefficient of~${\rm d}/{\rm d}t$ in~\eqref{eqn:eqproof1} is given by the rational function~$P'(t)/P(t)$, as in the statement~\eqref{eqn:f}.

The computation of the bracket $[f,f]_2$ needs a little more work. From the definition of RC brackets, we see that $[f,f]_2$ is a cusp form of weight eight. Moreover, it has a zero of order $2n-4$ where $f$ is zero, so it is necessarly divisible by~$f^2$. There exists then an element~$h_4\in M_4(\Gamma)$ such that $[f,f]_2=f^2h_4$. By Lemma~\ref{thm:dime} we know that~$h_4$ is of the form
\[
h_4=f^2Q(t),
\]
where $Q(t)$ is a polynomial in $t$ of degree $\dim{M_4(\Gamma)}=2n-3$. Since $[f,f]_2$ is a cusp form, $[f,f]_2/f^2$ has a zero in every cusp different from~$c_0$, and these zeros are simple. This means that the polynomial~$Q(t)$ is divisible by~$P(t)$. We have then
\[
h_4=f^2P(t)\bigl(\hat{\rho}_{n-3}t^{n-3}+\hat{\rho}_{n-2}t^{n-2}+\dots+\hat{\rho}_0\bigr),
\]
for some $\hat{\rho}_0,\dots,\hat{\rho}_{n-3}\in\C$. We can determine $\hat{\rho}_{n-3}$ by considering the expansion of $f$ at the cusp $c_0$. If $q_0$ denotes a local parameter at $c_0$, the expansions of~$f$ and~$t$ are given by
\[
f=cq_0^{n-2}+\cdots,\qquad t=sq_0^{-1}+\cdots,
\]
for some non-zero $c,s\in\C$.
In $c_0$ the bracket $[f,f]_2$ has expansion
\[
[f,f]_2=6ff''-9f'^2=-3c^2(n-2)^2q_0^{2n-4}+\cdots,
\]
while $h_4$ is given by
\[
h_4=\bigl(cq_0^{n-2}+\cdots\bigr)^2\bigl(\hat{\rho}_{n-3}s^{n-3}q_0^{3-n}+\cdots\bigr)\bigl(s^{n-1}q_0^{1-n}+\cdots\bigr)= \hat{\rho}_{n-3}c^2s^{2n-4}q_0^0+\cdots.
\]
The above expansions and the equality $[f,f]_2=f^2h_4$ imply that
\[
\hat{\rho}_{n-3}c^2s^{2n-4}=-(n-2)^2.
\]
From the relation $t'=\kappa^{-1}P(t)f$ we can compute the constant $\kappa$ in terms of the coefficients appearing in the expansions at $c_0$, obtaining $\kappa = -cs^{n-2}$. This implies that
\[
\hat{\rho}_{n-3}=-3\kappa^{-2}(n-2)^2.
\]
It finally follows that
\begin{align*}
-\frac{[f,f]_2}{12f^2t'^2}&=\frac{-f^4P(t)\big({-}3\kappa^{-2}(n-2)^2t^{n-3}+\hat{\rho}_4t^{n-4}+\cdots+\hat{\rho}_0\big)}{12\kappa^{-2}f^4P(t)^2}\\
 &=\frac{(n/2-1)^2t^{n-3}+\rho_{n-4}t^{n-4}+\cdots+\rho_0}{P(t)},
\end{align*}
where $\rho_i=-\hat{\rho_i}\kappa^2/12$.
\end{proof}

\subsection*{Acknowledgments}
The paper was written while I was a graduate student at SISSA (International School for Advanced Studies) in Trieste and a visiting student of the IMPRS graduate school at the Max Planck Institute for Mathematics in Bonn. I~want to thank both institutions for the excellent working conditions. I~want to thank my advisers Don Zagier and Fernando Rodriguez Villegas for their suggestions and support, and Giordano Cotti for very useful conversations. Finally I~thank the anonymous referees, whose suggestions and remarks gave a relevant contribution to improve the paper.

\pdfbookmark[1]{References}{ref}
\LastPageEnding


\begin{thebibliography}{99}
\footnotesize\itemsep=0pt

\bibitem{Anselmo}
Anselmo T., Nelson R., Carneiro~da Cunha B., Crowdy D.G., Accessory parameters
 in conformal mapping: exploiting the isomonodromic tau function for
 {P}ainlev\'e~{VI}, \href{https://doi.org/10.1098/rspa.2018.0080}{\textit{Proc.~A.}} \textbf{474} (2018), 20180080, 20~pages.

\bibitem{Beardon}
Beardon A.F., The geometry of discrete groups, \textit{Graduate Texts in
 Mathematics}, Vol.~91, \href{https://doi.org/10.1007/978-1-4612-1146-4}{Springer-Verlag}, New York, 1983.

\bibitem{Bers}
Bers L., Quasiconformal mappings, with applications to differential equations,
 function theory and topology, \href{https://doi.org/10.1090/S0002-9904-1977-14390-5}{\textit{Bull. Amer. Math. Soc.}} \textbf{83}
 (1977), 1083--1100.

\bibitem{Bogo}
Bogo G., Modular forms, deformation of punctured spheres, and extensions of
 symmetric tensor representations, \textit{Math. Res. Lett.}, {t}o appear,
 \href{https://arxiv.org/abs/#2}{arXiv:2004.04716}.

\bibitem{BouwMoeller}
Bouw I.I., M\"oller M., Differential equations associated with nonarithmetic
 {F}uchsian groups, \href{https://doi.org/10.1112/jlms/jdp059}{\textit{J.~Lond. Math. Soc.}} \textbf{81} (2010), 65--90,
 \href{https://arxiv.org/abs/0710.5277}{arXiv:0710.5277}.

\bibitem{123}
Bruinier J.H., van~der Geer G., Harder G., Zagier D., The 1-2-3 of modular
 forms, \href{https://doi.org/10.1007/978-3-540-74119-0}{Universitext}, Springer-Verlag, Berlin, 2008.

\bibitem{Chudnovsky}
Chudnovsky D.V., Chudnovsky G.V., Transcendental methods and theta-functions,
 in Theta Functions~-- {B}owdoin 1987, {P}art~2 ({B}runswick, {ME}, 1987),
 \textit{Proc. Sympos. Pure Math.}, Vol.~49, \href{https://doi.org/10.1103/physrevb.40.2946}{Amer. Math. Soc.}, Providence, RI,
 1989, 167--232.

\bibitem{SGervais}
de~Saint-Gervais H.P., Uniformization of {R}iemann surfaces, \textit{Heritage of
 European Mathematics}, \href{https://doi.org/10.4171/145}{European Mathematical Society (EMS)}, Z\"urich, 2016.

\bibitem{Ford}
Ford L.R., Automorphic functions, 2nd~ed., Chelsea Publishing Co., New York,
 1951.

\bibitem{Goldman}
Goldman W.M., Projective structures with {F}uchsian holonomy,
 \href{https://doi.org/10.4310/jdg/1214440978}{\textit{J.~Differential Geom.}} \textbf{25} (1987), 297--326.

\bibitem{Hempel}
Hempel J.A., On the uniformization of the {$n$}-punctured sphere, \href{https://doi.org/10.1112/blms/20.2.97}{\textit{Bull.
 London Math. Soc.}} \textbf{20} (1988), 97--115.

\bibitem{Hoffman}
Hoffman J., Monodromy calculations for some differential equations, Ph.D.~Thesis, {J}ohannes Gutenberg-Universit\"at Mainz, 2013.

\bibitem{KRV}
Keen L., Rauch H.E., Vasquez A.T., Moduli of punctured tori and the accessory
 parameter of {L}am\'e's equation, \href{https://doi.org/10.2307/1998172}{\textit{Trans. Amer. Math. Soc.}}
 \textbf{255} (1979), 201--230.

\bibitem{Kra}
Kra I., Accessory parameters for punctured spheres, \href{https://doi.org/10.2307/2001420}{\textit{Trans. Amer. Math.
 Soc.}} \textbf{313} (1989), 589--617.

\bibitem{Miyake}
Miyake T., Modular forms, \textit{Springer Monographs in Mathematics}, \href{https://doi.org/10.1007/3-540-29593-3}{Springer-Verlag},
 Berlin, 1989.

\bibitem{Nehari}
Nehari Z., On the accessory parameters of a {F}uchsian differential equation,
 \href{https://doi.org/10.2307/2372089}{\textit{Amer.~J. Math.}} \textbf{71} (1949), 24--39.

\bibitem{Poincare}
Poincar\'e H., Sur les groupes des \'equations lin\'eaires, \href{https://doi.org/10.1007/BF02418420}{\textit{Acta Math.}}
 \textbf{4} (1884), 201--312.

\bibitem{Polyakov}
Polyakov A.M., Quantum geometry of bosonic strings, \href{https://doi.org/10.1016/0370-2693(81)90743-7}{\textit{Phys. Lett.~B}}
 \textbf{103} (1981), 207--210.

\bibitem{Smirnov}
Smirnov V., Sur les \'equations diff\'erentielles lin\'eaires du second ordre et la th\'eorie des fonctions automorphes, \textit{Bull. Sci. Math.} \textbf{45} (1921), 93--120.

\bibitem{T2}
Takhtajan L.A., Liouville theory: quantum geometry of {R}iemann surfaces,
 \href{https://doi.org/10.1142/S0217732393002269}{\textit{Modern Phys. Lett.~A}} \textbf{8} (1993), 3529--3535,
 \href{https://arxiv.org/abs/hep-th/9308125}{arXiv:hep-th/9308125}.

\bibitem{Thompson}
Thompson J.G., Algebraic numbers associated to certain punctured spheres,
 \href{https://doi.org/10.1016/0021-8693(86)90234-6}{\textit{J.~Algebra}} \textbf{104} (1986), 61--73.

\bibitem{ZagApery}
Zagier D., Integral solutions of {A}p\'ery-like recurrence equations, in Groups
 and Symmetries, \textit{CRM Proc. Lecture Notes}, Vol.~47, \href{https://doi.org/10.1090/crmp/047/22}{Amer. Math. Soc.},
 Providence, RI, 2009, 349--366.

\bibitem{TZ1}
Zograf P.G., Takhtajan L.A., On {L}iouville equation, accessory parameters
 and the geometry of {T}eichm\"uller space for {R}iemann surfaces of
 genus~{$0$}, \href{https://doi.org/10.1070/SM1988v060n01ABEH003160}{\textit{Math. USSR-Sb.}} \textbf{60} (1988), 143--161.

\end{thebibliography}
\end{document}